\documentclass{article}
\usepackage{epsfig,latexsym,array,amsthm,amssymb,amsfonts,amsmath,multirow}
\usepackage[toc,page]{appendix}

\newtheorem{thm}{Theorem}[section]
\newtheorem{prop}[thm]{Proposition}
\newtheorem{lemma}[thm]{Lemma}
\newtheorem{cor}[thm]{Corollary}

\theoremstyle{definition}

\newtheorem{alg}[thm]{Algorithm}

\newcommand{\comment}[1]{} 

\newcommand{\lam}{\lambda}
\newcommand{\fs}{\left\lfloor\tfrac{s}{2}\right\rfloor}
\newcommand{\irange}{-\fs,\dotsc,(0),\dotsc,\fs}
\newcommand{\pf}{\noindent{\em Proof: }}
\def\qed{\hfill $\Box$ \vskip .15 in}
\makeatletter\def\imod#1{\allowbreak\mkern10mu({\operator@font mod}\,\,#1)}\makeatother

\numberwithin{equation}{section}

\date{\today}

\begin{document}

\title{\bf Nonexistence Results for Tight\\ Block Designs}

\author{
Peter Dukes and Jesse Short-Gershman\thanks{Research of the authors is supported by NSERC.}\\
Mathematics and Statistics\\
University of Victoria\\
Victoria, BC V8W 3R4\\
Canada \\
{\tt dukes@uvic.ca, jesseasg@uvic.ca}
}

\maketitle
\begin{abstract}
Recall that combinatorial $2s$-designs admit a classical lower bound $b \ge \binom{v}{s}$ on their number of blocks, and that a design meeting this bound is called tight.  A long-standing result of Bannai is that there exist only finitely many nontrivial tight $2s$-designs for each fixed $s \ge 5$, although no concrete understanding of `finitely many' is given.  Here, we use the Smith Bound on approximate polynomial zeros to quantify this asymptotic nonexistence.  Then, we outline and employ a computer search over the remaining parameter sets to establish (as expected) that there are in fact no such designs for $5 \le s \le 9$, although the same analysis could in principle be extended to larger $s$.  Additionally, we obtain strong necessary conditions for existence in the difficult case $s=4$.\\
\\
MSC Primary 05B05, 33D45; Secondary 05E30
\end{abstract}

\section{Introduction}

Let $v \ge k \ge t$ be positive integers and $\lam$ be a nonnegative integer.  A $t$-$(v,k,\lam)$ \emph{design}, or simply a $t$-\emph{design}, is a pair $(V,\mathcal B)$ where $V$ is a $v$-set and $\mathcal B$ is a collection of $k$-subsets of $V$ such that any $t$-subset of $V$ is contained in exactly $\lam$ elements of $\mathcal B$.  The elements of $V$ are \emph{points} and the elements of $\mathcal B$ are \emph{blocks}.  Since $t$-designs are also $i$-designs for $i \le t$, the parameter $t$ is typically called the \emph{strength}.  The number of blocks is usually denoted $b$ and an easy double-counting argument shows $b=\lam\binom{v}{t}/\binom{k}{t}$.  

Suppose $(V,\mathcal B)$ is a $t$-$(v,k,\lam)$ design. Generalizing Fisher's Inequality, Ray-Chaudhuri and Wilson \cite{RW} showed that if $t$ is even, say $t=2s$, and $v \ge k+s$, then $b \ge \binom{v}{s}$. If equality holds in this bound, we say $(V,\mathcal B)$ is \emph{tight}. The \emph{trivial} tight $2s$-designs are those with $v = k+s$, where each of the $\binom{v}{k} = \binom{v}{s}$ $k$-subsets of $V$ is a block.
The case of odd strength is investigated in \cite{IS}; however, it is impossible for $(2s-1)$-designs to be tight in the sense of having $\binom{v}{s-1}$ blocks.

Returning to even strength, the full set of parameters $v$ and $k$ for which a tight $2s$-design exists has only been determined for $s=2,3$. Note that, when $s=1$, tight 2-designs have $b=v$ and are the `symmetric' designs; see \cite{Ionin,Lander} for surveys of this rich (yet very challenging) topic.  In 1975, Ito \cite{Ito} published a proof that the only nontrivial tight 4-designs are the Witt 4-(23,7,1) design and its complementary 4-(23,16,52) design, but his proof was found to be incorrect. A few years later, Enomoto, Ito, and Noda \cite{EIN} proved the weaker result that there are finitely many nontrivial tight 4-designs, though still believing Ito's initial claim to be true.  Finally, in 1978, Bremner \cite{Bremner} successfully settled $s=2$ by reaffirming Ito's result.  Peterson \cite{Peterson} proved in 1976 that there exist no nontrivial tight 6-designs.

Bannai \cite{B} proved that there exist only finitely many nontrivial tight $2s$-designs for each $s \ge 5$. The case $s=4$ is quite open, and the `finitely many' for $s \ge 5$ is not explicit and potentially grows with $s$.  However, it is probably the case that there are no unknown tight $2s$-designs for $s \ge 2$.

Central to these negative results is the following strong condition, discovered first by Ray Chaudhuri and Wilson \cite{RW}, and also implicitly by Delsarte \cite{Delsarte}.

\begin{prop}
{\rm (\cite{Delsarte,RW})}
\label{psi}
If there exists a tight $2s$-$(v,k,\lam)$ design, then the zeros of the following degree $s$ polynomial $\Psi_s(x)$ are the intersection numbers of the design, and hence they must all be nonnegative integers:
\begin{equation}
\label{psi-def}
\Psi_s(x)=\sum_{i=0}^s(-1)^{s-i}\frac{\binom{v-s}{i}\binom{k-i}{s-i}\binom{k-1-i}{s-i}}{\binom{s}{i}}\binom{x}{i}.
\end{equation}
\end{prop}
The polynomials $\Psi_s$ are variants of the Hahn polynomials, \cite{Hahn}.  

Since a $2s$-design with $v \ge k+s$ induces at least $s$ intersection numbers \cite{RW}, it follows that the zeros of $\Psi_s$ must additionally be distinct integers for tight designs.  Note also that $\Psi_s$ has no dependence on $\lambda$; indeed, for tight designs $\lambda = \binom{v}{s}\binom{k}{2s} \binom{v}{2s}^{-1}$ and is therefore uniquely determined by $v$ and $k$.

Analogously, the Lloyd polynomials $L_e(x)$ are important for the characterization of perfect $e$-error-correcting codes; see \cite{vL}.  It is interesting that this characterization of perfect codes was completed long ago, while the open problems mentioned before Proposition~\ref{psi} remain for tight designs.  Our goal here is to revive the interest in tight designs and take a modest step toward the full characterization of their parameters. 

The outline is as follows.  In Section 2, we review the work of Bannai in \cite{B} on the asymptotic structure of the zeros of $\Psi_s$.  Extending this, we obtain some exact bounds relevant to this analysis.  Section 3 summarizes the techniques for exhausting small cases $s \ge 5$, and Section 4 is devoted to a partial analysis of the case $s=4$.  An appendix of tables following the main text will prove useful to the interested reader.

\section{Bannai's analysis and the Smith bound}

\subsection{Notation}

Assuming a tight design, let $x_i$, for $i = \irange$, denote the zeros of $\Psi_s$ listed in increasing order.  For example, the zeros of $\Psi_4$ and $\Psi_5$ are denoted $x_{-2}<x_{-1}<x_1<x_2$ and $x_{-2}<x_{-1}<x_0<x_1<x_2$, respectively.

An important parameter is the arithmetic mean of the zeros of $\Psi_s(x)$, which we denote by $\overline{\alpha}$.  
From the coefficient of $x^{s-1}$, we have
\begin{equation}
\label{alpha-bar}
\overline{\alpha} = \frac{(k-s+1)(k-s)}{v-2s+1}+\frac{s-1}{2}.
\end{equation}
Now define, as in \cite{B}, 
$$\alpha = \frac{(k-s+1)(k-s)}{v-2s+1},$$
so that $\overline{\alpha} = \alpha+(s-1)/2$.
Also, following Bannai's notation, let us redefine the parameter $t$ as
$$t= \frac{v-2s+1}{k-s+1}.$$
Note $t=2$ implies $v=2k+1$.  Moreover, if $v<2k$, we may complement blocks, replacing $k$ with $v-k$ and obtain $v>2k$.  This is discussed further in Section~\ref{t-sym}.

Finally, put $\beta = \left(1-\tfrac{1}{t}\right)\sqrt{\alpha}$.  In terms of $v$ and $k$,
$$\beta = \frac{(v-k-s)\sqrt{(k-s+1)(k-s)}}{(v-2s+1)^{3/2}}.$$
So, in particular, $\beta = 0$ if and only if $v=k+s$.  In some sense $\beta$ can be seen as measuring the `nontriviality' of a (tight) $2s$-design.
Note also that 
\begin{eqnarray}
\label{k}
k &=& t^3(t-1)^{-2} \beta^2 + s,~{\rm and}\\
\label{v}
v &=& t^4(t-1)^{-2} \beta^2 +t+2s-1.
\end{eqnarray}

Bannai's proof of the existence of only finitely many nontrivial tight $2s$-designs, $s \ge 5$, is divided into cases according to this parameter $\beta$.  In particular, he proves
\begin{itemize}
\item
for any $\beta_0$, there are only finitely many tight $2s$-designs with $\beta \le \beta_0$; and
\item
there exists $\beta_0$ (depending only on $s$), such that there are no nontrivial tight $2s$-designs with $\beta > \beta_0$.
\end{itemize}

Here, our main goal is to compute such a $\beta_0$ explicitly for $5 \le s \le 9$ and, by searching across all pairs $(v,k)$ for which $\beta \le \beta_0$, show that there are in fact zero nontrivial tight $2s$-designs for these $s$.

\subsection{Symmetry with respect to the parameter $t$}
\label{t-sym}

In the analytic work which follows, it is helpful to obtain a lower bound on $t$.  As discussed above, we may complement blocks to assume $v \ge 2k$. The following was mentioned but not fully proven in \cite{B}.

\begin{lemma}
\label{vnot2k}
Let $s \ge 1$. There are no tight $2s$-designs with $v = 2k$.
\end{lemma}
\pf
Suppose $v = 2k$. Then from (\ref{alpha-bar}), $s\overline{\alpha}-\binom{s}{2} = s\alpha = \frac{s(k-s+1)(k-s)}{2(k-s)+1}$. Without too much effort, it can be seen that the least residue of $s(k-s+1)(k-s) \pmod{2(k-s)+1}$, denoted here by $r_{k,s}$, satisfies
$$r_{k,s}=
\begin{cases}
2(k-s)-\frac{s-4}{4} & \text{if $s \equiv 0 \imod{4}$};\\
k-s-\frac{s-2}{4} & \text{if $s \equiv 2 \imod{4}$};\\
\frac{k-s}{2}-\frac{s-1}{4} & \text{if $s \equiv 1 \imod{4}$ and $k$ is odd} \\
 & \text{~~~or $s \equiv 3 \imod{4}$ and $k$ is even};\\
\frac{3}{2}(k-s)-\frac{s-3}{4} & \text{if $s \equiv 3 \imod{4}$ and $k$ is odd}\\
 & \text{~~~or $s \equiv 1 \imod{4}$ and $k$ is even}.\\
\end{cases}$$
Since $k-s \ge s$, it follows that in all cases $r_{k,s}$ is an integer lying strictly between 0 and $2(k-s)+1$, so $s\overline{\alpha}-\binom{s}{2}$ is not an integer. But the integrality of $s\overline{\alpha}$ is necessary for the existence of a tight design since it is the sum of the zeros of $\Psi_s(x)$; therefore there are no tight $2s$-designs with $v = 2k$.
\qed

Now, we are able to justify assuming that $t \ge 2$ for nonexistence of tight designs.

\begin{prop}
\label{t}
Let $s \ge 1$. If there exists a nontrivial tight $2s$-design with $t < 2$, then there also exists a nontrivial tight $2s$-design with $t \ge 2$.
\end{prop}
\pf
Suppose $\mathcal D$ is a nontrivial tight $2s$-$(v,k,\lam)$ design with $t < 2$. This means $k \le v \le 2k-1$ because $v \neq 2k$ by Lemma~\ref{vnot2k}, and so the complementary $2s$-$(v,v-k,\lam')$ design of $\mathcal D$ is a nontrivial tight $2s$-design with $t \ge 2$.
\qed

Incidentally, Bannai and Peterson ruled out the case $t=2$, observing that it yields symmetric zeros of $\Psi_s$ about their mean $\overline{\alpha}$.  This is a key observation.

\begin{prop}
\label{v2k+1}
{\rm (\cite{B,Peterson})}
There does not exist any tight $2s$-design with $v=2k+1$.
\end{prop}

\subsection{Hermite polynomials}

Let $H_s(x)$ denote the normalized Hermite polynomial of degree $s$ defined recursively by $H_0(x) = 1$, $H_1(x) = x$, and for $s \ge 2$,
\begin{equation*}
H_s(x) = xH_{s-1}(x)-(s-1)H_{s-2}(x).
\end{equation*} 
Furthermore, let $\xi_i$, $i = \irange$, denote the zeros of $H_s(x)$ listed in increasing order. It is easily seen that $\xi_{-i} = -\xi_i$ for each $i$. See Appendix~\ref{hzeros} for a table of zeros of $H_s(x)$, $1 \le s \le 10$.  
In particular, for the analytical work in Section 3, we will make use of the following known estimates.

\begin{prop}
\label{xi}~
\begin{enumerate}
\item[{\rm (i)}] If $s$ is odd and $\ge 5$, then $\xi_1^2 < \sqrt{3}$.
\item[{\rm (ii)}] If $s$ is even and $\ge 8$, then $\xi_2^2-\xi_1^2 < \sqrt{3}$.
\item[{\rm (iii)}] If $s=6$, then $1.0 < \frac{\xi_2^2-\xi_1^2}{3} < 1.1$, $3.5 < \frac{\xi_3^2-\xi_1^2}{3} < 3.6$, and $3.34634 < \frac{\xi_3^2-\xi_1^2}{\xi_2^2-\xi_1^2} < 3.34635$.
\end{enumerate}
\end{prop}
\pf
Items (i) and (ii) are referenced in Bannai's Proposition 13 and proven on page 126 of \cite{Szego}. Item (iii) can be verified numerically. See Appendix~\ref{hzeros}. (Note that Bannai's Proposition 13 (iii) actually contains an error).
\qed

A useful identity is
\begin{equation}
\label{Hermite-derivative}
H_s'(x) = sH_{s-1}(x).
\end{equation}

For later reference we define, again as in \cite{B}, 
\begin{equation}
\label{lam-i-def}
\lam_i=\lam_i(t) = \left(1-\tfrac{2}{t} \right)^2 \left(\tfrac{\xi_i^2}{6} - \tfrac{s-1}{6}\right).
\end{equation}
Informally, Proposition 16 in \cite{B} states that as $\beta\to\infty$, the zeros $x_i$ of $\Psi_s(x)$ approach $\overline{\alpha}+\beta\xi_i+\lam_i$.  That is, when suitably normalized, $\Psi_s$ behaves like $H_s$ for large $\beta$ and fixed $t$.

\subsection{The Smith bound}

We now state a useful result for explicitly finding $\beta_0$.  Sometimes known as the Smith bound, it is a consequence of the Gershgorin circle theorem.

\begin{thm}{\rm (\cite{Smith})}
\label{smithbound}
Let $P(z)$ be a monic polynomial of degree $n$ and let $\xi_1,\dots,\xi_n$ be distinct points approximating the zeros of $P(z)$.  Define the circles
\[\Gamma_i = \left\{z : |z - \xi_i| \le \frac{n|P(\xi_i)|}{|Q'(\xi_i)|}\right\},\]
where $Q(z)$ is the monic polynomial of degree $n$ with zeros $\xi_1,\dots,\xi_n$.  Then the union of the circular regions $\Gamma_i$ contains all the zeros of $P(z)$, and any connected component consisting of just $k$ circles $\Gamma_i$ contains exactly $k$ zeros of $P(z)$.
\end{thm}

Let $s \ge 1$.  For each $i \in \{\irange\}$, define the monic degree $s$ polynomial (in $z$) 
\begin{equation}
\label{G-def}
G_s^{(i)}(z) = \frac{s!}{\beta^s \binom{v-s}{s}} \Psi_s(\overline{\alpha} + \beta z + \lam_i),
\end{equation}
and put $z_i = (x_i-\overline{\alpha}-\lam_i)/\beta$, the zero of $G_s^{(i)}(z)$ corresponding to $x_i$.

We will see from Propositions \ref{smithtight} and \ref{beta1} that the $z_i$ are well-approximated by the $\xi_i$ as $\beta\to\infty$, independently of $t$. 

\begin{prop}
\label{smithtight}
Let $s \ge 1.$  Then
\[|z_i - \xi_i| \le \frac{|G_s^{(i)}(\xi_i)|}{|H_{s-1}(\xi_i)|}.\]
\end{prop}
\pf
Simply apply Theorem~\ref{smithbound} to the polynomial $G_s^{(i)}(z)$, letting $Q(z) = H_s(z)$, to get
\[|z_i - \xi_i| \le \frac{ s |G_s^{(i)}(\xi_i)| }{ |H_s'(\xi_i)| }.\]   The result now follows from (\ref{Hermite-derivative}).
\qed

\subsection{Bounding $G_s$ in terms of $\beta$}

In the next proposition, it is helpful to think of the $G_s^{(i)}(\xi_i)$ as functions of $\beta$ and $t$. 

\begin{prop}
\label{beta1}
Let $s \ge 2$. For each $i \in \{\irange\}$, there exist constants $B_i,C_i$ such that whenever $\beta > B_i$,
$$|G_s^{(i)}(\xi_i)| < \frac{C_i}{\beta^2}$$
for all $t \ge 2$.
\end{prop}

The necessary ingredients for this result were proved in \cite{B}, although the bound was not directly stated in this form.  Therefore, we omit the proof and instead focus on how to (carefully) obtain $B_i$ and $C_i$ for small $s$ using some basic computer algebra.

\begin{alg}  
\label{al}
For fixed $s$ and $i$, we may obtain constants $B_i$ and $C_i$ in Proposition~\ref{beta1} by the following procedure.
\begin{enumerate}
\item
Using (\ref{psi-def}), substitute (\ref{alpha-bar}), (\ref{k}), (\ref{v}) and (\ref{lam-i-def}) into (\ref{G-def}).  To defer floating-point precision issues, we first replace $\xi_i$ in (\ref{lam-i-def}) by a symbolic parameter $r$.
\item
This results in an expression for $G_s^{(i)}(r)$ as a rational function of $\beta$, say
$$G_s^{(i)}(r)(\beta,t) = \frac{p(r,\beta,t)}{q(\beta,t)}.$$
Here, the denominator is 
\begin{equation}
\label{den-def}
q(\beta,t) = \beta^s \tbinom{v-s}{s} = \beta^s \binom{t^4(t-1)^{-2} \beta^2 +t+s-1}{s}.
\end{equation}
\item
Observe that $q$ is positive for $\beta>0$ and $t \ge 2$, and that a lower bound on $q$ is
$$\tilde{q}(\beta,t) = \tfrac{1}{s!} \beta^{3s} t^{4s} (t-1)^{-2s}.$$
This is obtained by replacing each factor in the falling factorial of (\ref{den-def}) by $t^4(t-1)^{-2} \beta^2$.
\item
The numerator $p(r,\beta,t)$ is, for general $r$, a polynomial of degree $3s$ in $\beta$.  However, for $r=\xi_i$, Proposition~\ref{beta1} shows the two top coefficients, namely of $\beta^{3s}$ and $\beta^{3s-1}$, vanish.  Again, to maintain symbolic algebra, we artifically replace these coefficients by zero and call this polynomial $\tilde{p}(r,\beta,t)$. 
\item
We have 
$$\beta^2 G_s^{(i)}(r)(\beta,t) \le \frac{\beta^2 \tilde{p}(r,\beta,t)}{\tilde{q}(\beta,t)}.$$
Note that for $r=\xi_i$, the right hand side is a polynomial in $\beta^{-1}$.  
\item
Consider the coefficient $\kappa_j(r,t)$ of $\beta^{3s-j}$ in $\beta^2 \tilde{p}(r,\beta,t)$.
With $r=\xi_i$, compute (or upper-bound) the maxima
$$M_{j}= \sup_{t \ge 2} \frac{|\kappa_j(\xi_i,t)|}{\tfrac{1}{s!}(t-1)^{-2s} t^{4s}}.$$
Then, estimating term-by-term,
$$|\beta^2 G_s^{(i)}(\xi_i)(\beta,t)| \le M_0+M_1 \beta^{-1} + M_2 \beta^{-2} + \dots$$
for all $t \ge 2$.  
\item
Construct $B_i,C_i$ so that $\beta > B_i$ implies $M_0+ M_1\beta^{-1} + M_2 \beta^{-2} + \dots \le C_i$.  Note that with sufficiently large $B_i$ and a safe choice of $C_i$, it suffices to estimate the first few coefficients $M_j$.  
\end{enumerate}
\end{alg}

We should remark that for small $s$, Algorithm~\ref{al} -- even the calculation of all $3s-1$ coefficient maxima $M_j$ -- is essentially instantaneous on today's personal computers.  Moreover, deferring the use of floating-point arithmetic to step 5 -- when $t$ is eliminated -- makes our subsequent use of floating-point numbers $M_j$ quite mild.  Indeed, there is virtually no loss in taking $M_j$ as (integer) ceilings of the suprema, so that estimating for $C_i$ can be performed in $\mathbb{Q}$.

See Appendix~\ref{C_is} for the results of this calculation for each $5 \le s \le 9$ and all relevant indices $i$.

\subsection{Bounding the zeros}

We are now ready for our main result of this section.  This is in Bannai's paper \cite{B}, but with no attempt to control $\beta$.

\begin{prop}
\label{beta2}
Fix a positive integer $s$ and $i \in \{\irange\}$.  Put 
$y_i = x_i - \overline{\alpha} - \beta \xi_i,$ where recall $x_i$ and $\xi_i$ are corresponding roots of $\Psi_s$ and $H_s$, respectively.  Let $\epsilon > 0$ and define 
$$\widehat{\beta}(i,\epsilon) = \max\left\{B_i,\frac{C_i}{\epsilon D_i}\right\},$$
where $D_i=|H_{s-1}(\xi_i)|$.  Then for all $\beta > \widehat{\beta}$ and all $t \ge 2$, 
$$|y_i - \lam_i| < \epsilon.$$
\end{prop}
\pf
Observe that $|y_i-\lam_i| = \beta|z_i-\xi_i|$, since
$$x_i=\overline{\alpha}+\beta\xi_i+y_i=\overline{\alpha}+\beta z_i+\lam_i.$$  
The estimate now follows easily from Propositions~\ref{smithtight} and \ref{beta1}.
\qed

\section{The case $s \ge 5$}

\subsection{Estimates for large $\beta$}

The goal here is to provide formulas for the smallest $\beta_0$ possible (see the end of Section 2.1) using the $B_i$ and $C_i$ constructed in Algorithm \ref{al}. This task is simplified under the conditions that $B_i$ is independent of $i$ and $C_i=C_{-i}$. There is no loss of generality in assuming this because we can simply take $\overline{B}$ to be the maximum of the $B_i$ and $\overline{C_i}=\max\{C_i,C_{-i}\}$, and then redefine each $B_i = \overline{B}$ and $C_i = C_{-i} = \overline{C_i}$. In fact, this is not necessary for our explicit constructions because the constants in Appendix~\ref{C_is} satisfy the above conditions.

Again, for convenience, we denote $|H_{s-1}(\xi_i)|$ by $D_i$ in the following proofs.

\begin{prop}
\label{sOdd}
Let $s \ge 5$ be odd.
\begin{enumerate}
\item[{\rm (i)}] There exists $\beta_1$ such that, whenever $\beta > \beta_1$,
\[|y_1+y_{-1}-2y_0| < 1.\]
\item[{\rm (ii)}] There exists $\beta_2$ such that, whenever $\beta > \beta_2$,
\[|y_i+y_{-i}-y_{i-1}-y_{-(i-1)}| < 1+\frac{\xi_i^2-\xi_{i-1}^2}{\xi_{i-1}^2-\xi_{i-2}^2}|y_{i-1}+y_{-(i-1)}-y_{i-2}-y_{-(i-2)}|\]
for $2\le i\le\fs$.
\item[{\rm (iii)}] There exists $\beta_0(s)$ such that, whenever $\beta > \beta_0$ and $y_i+y_{-i}-y_{i-1}-y_{-(i-1)}$ is an integer for $1\le i\le\fs$, it is necessarily the case that $y_i+y_{-i}-y_{i-1}-y_{-(i-1)} = 0$ for each $i$.
\end{enumerate}
\end{prop}
\pf
\begin{enumerate}
\item[{\rm (i)}] Observe
%$$y_1+y_{-1}-2y_0 = (y_1-\lam_1)+(y_{-1}-\lam_{-1})-2(y_0-\lam_0)+2(\lam_1-\lam_0),$$
that since $t \ge 2$, 
\begin{equation}
\label{sOdd-est}
0 \le 2(\lam_1-\lam_0) = (1-\tfrac{2}{t})^2\tfrac{\xi_1^2}{3} <\tfrac{\xi_1^2}{3}.
\end{equation}
Define
\[\epsilon_0 = \frac{1}{2}\left(1-\frac{\xi_1^2}{3}\right)\left(1+\frac{C_1D_0}{C_0D_1}\right)^{-1}\hspace{1mm}\text{and}\hspace{3mm}\beta_1 = \widehat{\beta}(0,\epsilon_0).\]
If $\epsilon_1=\epsilon_0\tfrac{C_1D_0}{C_0D_1}$, then 
\[\widehat{\beta}(1,\epsilon_1)=\beta_1\hspace{3mm}\text{and}\hspace{3mm}2\epsilon_0+2\epsilon_1=1-\frac{\xi_1^2}{3}.\]
%%latestmod  -- both of the above statements are what imply the next line.  There is no point in numbering the equation in my opinion.
Hence for $\beta > \beta_1$,
\begin{eqnarray*}
|y_1+y_{-1}-2y_0-2(\lam_1-\lam_0)| &\le& |y_1-\lam_1|+|y_{-1}-\lam_{-1}|+2|y_0-\lam_0|, \\
& < & 2\epsilon_0+2\epsilon_1 = 1-\frac{\xi_1^2}{3},
\end{eqnarray*}
By (\ref{sOdd-est}), 
$$-\left(1-\tfrac{\xi_1^2}{3}\right) < y_1+y_{-1}-2y_0 < 1$$
and the claim follows.

\item[{\rm (ii)}] For $2\le i\le\fs$, let 
\[a_i=\frac{\xi_i^2-\xi_{i-1}^2}{\xi_{i-1}^2-\xi_{i-2}^2}\hspace{3mm}\text{and}\hspace{3mm}\epsilon_i = \frac{1}{2}\left(1+(1+a_i)\frac{C_{i-1}D_i}{C_iD_{i-1}}+a_i\frac{C_{i-2}D_i}{C_iD_{i-2}}\right)^{-1}.\]
Note if $2\le i\le\fs$, then
\begin{equation}
\label{lam-rel}
(\lam_i-\lam_{i-1}) = (\lam_{i-1}-\lam_{i-2})a_i.
\end{equation}
Define $\beta_2 = \max\left\{\widehat{\beta}(i,\epsilon_i):2\le i\le\fs\right\}$. 
For $\beta > \beta_2$ and working as in (i),
\begin{equation}
\label{xx}
|y_i+y_{-i}-y_{i-1}-y_{-(i-1)}-2(\lam_i-\lam_{i-1})| < 2\epsilon_i \left( 1+\tfrac{C_{i-1}D_i}{C_iD_{i-1}} \right).
\end{equation}
Using (\ref{lam-rel}) and (\ref{xx}) again with $i-1$ replacing $i$,
\begin{equation}
\nonumber
\begin{split}
|y_i+y_{-i}-&y_{i-1}-y_{-(i-1)}| \\
 &< 2\epsilon_i \left(1+\tfrac{C_{i-1}D_i}{C_iD_{i-1}}\right) + 2 \epsilon_i a_i \left(\tfrac{C_{i-1}D_i}{C_iD_{i-1}} + \tfrac{C_{i-2}D_i}{C_iD_{i-2}} \right) \\
 &~~~~~~~~~~ + a_i|y_{i-1}+y_{-(i-1)}-y_{i-2}-y_{-(i-2)}| \\
 &=1+a_i|y_{i-1}+y_{-(i-1)}-y_{i-2}-y_{-(i-2)}|,
\end{split}
\end{equation}
as required.

\item[{\rm (iii)}] Set $\beta_0(s) = \max\{\beta_1,\beta_2\}$ and assume that $\beta > \beta_0(s)$ and $y_i+y_{-i}-y_{i-1}-y_{-(i-1)}$ is an integer for $1\le i\le\fs$. By (i), $y_1+y_{-1}-2y_0 = 0$ since it is an integer whose absolute value is less than 1. Assume that $y_{i-1}+y_{-(i-1)}-y_{i-2}-y_{-(i-2)} = 0$ for some $2\le i\le\fs$. Then (ii) gives that $|y_i+y_{-i}-y_{i-1}-y_{-(i-1)}|$ is also less than one and hence equal to 0 since it is an integer, so by induction $y_i+y_{-i}-y_{i-1}-y_{-(i-1)} = 0$ for $1\le i\le\fs$, and so the proof is complete. \qed
\end{enumerate}

\begin{prop}
\label{sEven}
Let $s \ge 8$ be even.
\begin{enumerate}
\item[{\rm (i)}] There exists $\beta_1$ such that, whenever $\beta > \beta_1$,
\[|y_2+y_{-2}-y_1-y_{-1}| < 1.\]
\item[{\rm (ii)}] There exists $\beta_2$ such that, whenever $\beta > \beta_2$,
\[|y_i+y_{-i}-y_{i-1}-y_{-(i-1)}| < 1+\frac{\xi_i^2-\xi_{i-1}^2}{\xi_{i-1}^2-\xi_{i-2}^2}|y_{i-1}+y_{-(i-1)}-y_{i-2}-y_{-(i-2)}|\]
for $3\le i\le\fs$.
\item[{\rm (iii)}] There exists $\beta_0(s)$ such that, whenever $\beta > \beta_0(s)$ and $y_i+y_{-i}-y_{i-1}-y_{-(i-1)}$ is an integer for $2\le i\le\fs$, it is necessarily the case that $y_i+y_{-i}-y_{i-1}-y_{-(i-1)} = 0$ for each $i$.
\end{enumerate}
\end{prop}
\pf
\begin{enumerate}

\item[{\rm (i)}] Since $t \ge 2$, 
\begin{equation}
\label{sEven-est}
0 \le 2(\lam_2-\lam_1) = (1-\tfrac{2}{t})^2\tfrac{\xi_2^2-\xi_1^2}{3} < \tfrac{\xi_2^2-\xi_1^2}{3}.
\end{equation}
Define
\[\epsilon_1 = \frac{1}{2}\left(1-\frac{\xi_2^2-\xi_1^2}{3}\right)\left(1+\frac{C_2D_1}{C_1D_2}\right)^{-1}\hspace{1mm}\text{and}\hspace{3mm}\beta_1 = \widehat{\beta}(1,\epsilon_1).\]
If $\epsilon_2=\epsilon_1\tfrac{C_2D_1}{C_1D_2}$, then
\[\widehat{\beta}(2,\epsilon_2)=\beta_1\hspace{3mm}\text{and}\hspace{3mm}2\epsilon_1+2\epsilon_2=1-\frac{\xi_2^2-\xi_1^2}{3}.\]
%%latestmod  -- see the comments in Prop \ref{sOdd} (i)
Hence for $\beta > \beta_1$,
\[|y_2+y_{-2}-y_1-y_{-1}-2(\lam_2-\lam_1)| < 2\epsilon_1+2\epsilon_2=1-\frac{\xi_2^2-\xi_1^2}{3}.\]
By (\ref{sEven-est}), 
$$-\left(1-\tfrac{\xi_2^2-\xi_1^2}{3}\right) < y_2+y_{-2}-y_1-y_{-1} < 1$$
and the claim follows.

\item[{\rm (ii)}] Define $\epsilon_i$ and $\beta_2$ in as in the proof of Proposition \ref{sOdd} (ii), but omit $i=2$.
\item[{\rm (iii)}] Imitate the proof of Proposition \ref{sOdd} (iii). \qed
\end{enumerate}

In the case $s=6$, $\tfrac{\xi_2^2-\xi_1^2}{3} > 1$. Hence it is impossible to choose a $\beta_1$ to guarantee that $y_2+y_{-2}-y_1-y_{-1} = 0$ whenever it is an integer and $\beta > \beta_1$.%%latestmod

\begin{prop}
\label{sEquals6}
Let $s=6$. There exists $\beta_0(6)$ such that, whenever $\beta > \beta_0(6)$ and $(y_2+y_{-2}-y_1-y_{-1})$, $(y_3+y_{-3}-y_1-y_{-1})$ are both integers, it is necessarily the case that $y_2+y_{-2}-y_1-y_{-1} = y_3+y_{-3}-y_1-y_{-1} = 0$.
\end{prop}
\pf
Observe
\[0 \le 2(\lam_2-\lam_1) < \tfrac{\xi_2^2-\xi_1^2}{3} < 1.1~~\text{ and}~~
0 \le 2(\lam_3-\lam_1) < \tfrac{\xi_3^2-\xi_1^2}{3} < 3.6\]
by Proposition~\ref{xi} (iii). Let $a = \tfrac{\xi_3^2-\xi_1^2}{\xi_2^2-\xi_1^2}$ and define
\[\epsilon_1 = \frac{1}{2}\left(a-3\right)\left(1+a+a\frac{C_2D_1}{C_1D_2}+\frac{C_3D_1}{C_1D_3}\right)^{-1}\hspace{1mm}\text{and}\hspace{3mm}\beta_0(6) = \widehat{\beta}(1,\epsilon_1).\]
Then, with
\[\epsilon_2 = 2\epsilon_1\left(1+\frac{C_2D_1}{C_1D_2}\right)\hspace{3mm}\text{and}\hspace{3mm}\epsilon_3 = 2\epsilon_1\left(1+\frac{C_3D_1}{C_1D_3}\right),\]
we have
\[\epsilon_2a+\epsilon_3 = a-3,\]
\[0 < \epsilon_2 < (a-3)/a\approx 0.10350,\hspace{3mm}\text{and}\]
\[0 < \epsilon_3 < a-3\approx 0.34635.\]
Assume $\beta > \beta_0(6)$. Then
$|y_2+y_{-2}-y_1-y_{-1}-2(\lam_2-\lam_1)| < \epsilon_2$ implies 
$$y_2+y_{-2}-y_1-y_{-1} \in \{0,1\}.$$
Likewise,
$|y_3+y_{-3}-y_1-y_{-1}-2(\lam_3-\lam_1)| < \epsilon_3$ implies
\begin{equation}
\label{0123}
y_3+y_{-3}-y_1-y_{-1} \in \{0,1,2,3\}.
\end{equation}
If $y_2+y_{-2}-y_1-y_{-1} = 0$, then $2(\lam_2-\lam_1) < \epsilon_2$ and so $0 \le 2(\lam_3-\lam_1) < \epsilon_2a$. Hence $y_3+y_{-3}-y_1-y_{-1} < \epsilon_2a+\epsilon_3 = a-3\approx 0.34635$, and so $y_3+y_{-3}-y_1-y_{-1} = 0$.
On the other hand, suppose $y_2+y_{-2}-y_1-y_{-1} = 1$. Then $2(\lam_2-\lam_1) > 1-\epsilon_2$ and so $2(\lam_3-\lam_1) > (1-\epsilon_2)a = a-\epsilon_2a$. Hence $y_3+y_{-3}-y_1-y_{-1} > a-\epsilon_2a-\epsilon_3 = 3$, a contradiction to (\ref{0123}).

\medskip
\noindent
It follows that $y_2+y_{-2}-y_1-y_{-1} = y_3+y_{-3}-y_1-y_{-1} = 0$.
\qed

To summarize, we have the following reworking of Proposition 17 in \cite{B}, but with explicit $\beta_0$.

\begin{thm}
\label{bannai}
For each $s \ge 5$, there are no tight $2s$-designs with $\beta > \beta_0(s)$.
\end{thm}
\pf
Suppose $x_{-\fs} < \cdots < x_{\fs}$  are the intersection numbers of a tight $2s$-design with $\beta > \beta_0(s)$. By Proposition~\ref{t}, we may assume $t \ge 2$. Then, since $\xi_{-i}=-\xi_i$ and $\lam_{-i}=\lam_i$, we have $x_i+x_{-i}-x_j-x_{-j} = y_i+y_{-i}-y_j-y_{-j}$, and this implies that $y_i+y_{-i}-y_j-y_{-j}$ is an integer for each $i,j \in \{(0),1,2,\dots,\fs\}$.  By Propositions \ref{sOdd} (iii), \ref{sEven} (iii) and \ref{sEquals6}, these integers must vanish.  Specifically,

\medskip
\noindent
\underline{\emph{Case 1}:} $s$ is odd and $\ge 5$ implies $y_i+y_{-i}-y_{i-1}-y_{-(i-1)} = 0$ for $1\le i\le\fs$.\\
\underline{\emph{Case 2}:} $s$ is even and $\ge 8$ implies $y_i+y_{-i}-y_{i-1}-y_{-(i-1)} = 0$ for $2\le i\le\fs$.\\
\underline{\emph{Case 3}:} $s=6$ implies $y_2+y_{-2}-y_1-y_{-1} = y_3+y_{-3}-y_1-y_{-1} = 0$.

\medskip
\noindent
In each case, the $x_i$ are symmetric about their arithmetic mean $\overline{\alpha}$. By Proposition 2 in \cite{B}, this implies $v=2k+1$.  Proposition~\ref{v2k+1} says this is impossible, and the proof is therefore complete.  
\qed

\subsection{Searching over small $\beta$}

We now turn to small values of $\beta$, for which the problem becomes finite.

\begin{alg}
\label{search}
To exclude tight $2s$-designs with $\beta \le \beta_0$, we may implement the following steps.
\begin{enumerate}
\item
Compute $\beta_0$ from the $B_i,C_i$ as in the previous section.
\item
By Propositions~\ref{t} and \ref{v2k+1}, we may restrict attention to $t>2$.  Since $\alpha = \beta^2/(1-\frac{1}{t})^2$, it follows that $\alpha < 4\beta_0^2$.  Now, since $\alpha = \left(s\overline{\alpha}+\binom{s}{2}\right)/s$ and $s\overline{\alpha}$ is an integer, we have $\alpha \in \tfrac{1}{s} \mathbb{Z}$.  This gives an explicit finite number of admissible $\alpha$, as Bannai observed in \cite{B}.%%latestmod
\item
Note that, under the assumption of a tight design, the expression
\begin{equation}
\label{coef-integer}
\binom{s}{2}\alpha\left(\alpha+\frac{2\alpha t-\alpha+2}{\alpha t^2+t+1}\right)
\end{equation}
is an integer.  This is because Proposition 5 in \cite{B} asserts that the coefficient of $x^{s-2}$ in the monic polynomial $s!\Psi_s(x)/\binom{v-s}{s}$ is
$$\binom{s}{2}\alpha\left(\alpha+\frac{2\alpha t-\alpha+2}{\alpha t^2+t+1}\right)+\binom{s}{3}\left(3\alpha+\frac{3s-1}{4}\right),$$
and the latter term is always an integer.
\item
Fix $\alpha$ as in Step 2.  Put $n = k-s = \alpha t$ and define 
$$g_{\alpha}(n) := \binom{s}{2}\alpha^2\left(1+\frac{2n-\alpha+2}{n^2+n+\alpha}\right)$$
as in (\ref{coef-integer}).
As $t>2$, we may take a lower bound $n_{\min}(\alpha) = \max\{s,\lfloor 2\alpha \rfloor+1\}$.  
\item
Since $g'_{\alpha}(n) < 0$ for all $n \ge n_{\min}(\alpha)$, it suffices to loop on integers $n$ from $n_{\min}(\alpha)$ until $n_{\max}(\alpha)$, where 
$g_\alpha(n_{\max}(\alpha)) \le \lfloor \binom{s}{2}\alpha^2 \rfloor+1.$
Any pairs $(k,v)$ which give integral $g_\alpha(n)$ are obtained by $k = n+s$ and $v = \frac{n^2+n}{\alpha}+2s-1$.
\item
In principle, at this point the zeros of $\Psi_s$ for these pairs $(k,v)$ can be analyzed.  However, in practice we found it sufficient in all cases to merely see that $\lam = \binom{v}{s}\binom{k}{2s}/\binom{v}{2s}$ was never even an integer.
\end{enumerate}
\end{alg}
 
We wrote a C program that implements Algorithm~\ref{search} for a given $s$ and $\beta_0$, but with an important optimization. For $n$ near $n_{\max}(\alpha)$, $|g'_{\alpha}(n)|$ is very small so it would be inefficient to loop over $n$ in this region. Therefore, the program loops over integer values of $g_\alpha(n)$ from $\lceil g_\alpha(n_{\max}(\alpha))\rceil$ and checks the integrality of the corresponding $n$ until the derivative becomes larger than a certain threshold (in absolute value), at which point it begins looping over $n$ to a much smaller $n_{\max}$. The program is available by contacting the authors.

Our calculations of $\beta_0$ in step 1 are displayed in Appendix~\ref{C_is}.  We can report that the method succeeds for $5 \le s \le 9$, and probably much higher $s$.  We have chosen to avoid continued searches for $s>9$ until new ideas are obtained.  In particular, it would be interesting if $s \ge s_0$ could be excluded for nontrivial tight $2s$-designs.

\begin{thm}
\label{none}
For each $5 \le s \le 9$, there are no nontrivial tight $2s$-designs.
\end{thm}

\section{The case $s = 4$}

The same analytic approach that is successful for $s \ge 5$ fails when $s=4$.
We can only guarantee that 
$$0 \le 2(\lam_2-\lam_1) < \frac{\xi_2^2-\xi_1^2}{3} = \sqrt{8/3} \approx 1.63299$$ when $t \ge 2$, and so there does not exist $\beta_0$ such that $|y_2+y_{-2}-y_1-y_{-1}| < 1$ for all $\beta > \beta_0$. 

However, it is possible to bound $|y_2+y_{-2}-y_1-y_{-1}|$ away from 2. Let
\[\epsilon_1=\frac{1}{2}\left(2-\sqrt{8/3}\right)\left(1+\frac{C_2D_1}{C_1D_2}\right)^{-1}\hspace{1mm}\text{and}\hspace{3mm}\beta_\star(4)=\widehat{\beta}(1,\epsilon_1).\]
%%latestmod
Then the existence of a tight 8-design with $\beta > \beta_\star(4)$ and $t \ge 2$ implies $y_2+y_{-2}-y_1-y_{-1} = 1$ and $2(\lam_2-\lam_1) \approx 1$, for which 
$$t = \tfrac{v-7}{k-3} \approx \frac{2}{1-\sqrt[4]{3/8}} \approx 9.1971905725.$$ 

We are able to obtain more precise conditions in the following result.

\begin{prop}
\label{sEquals4}
If there exists a nontrivial tight $8$-design with parameters $v$ and $k$, then $k >$ {\rm 25,000} and $f(k,v) = 0$, where $f(k,v)$ is as in Appendix~\ref{fgkv}.
\end{prop}
\pf
We first used Algorithm~\ref{search} to find that there are no nontrivial tight 8-designs with $\beta \le \beta_\star(4)$.  Thus, any tight 8-design with $t \ge 2$ must have $x_2+x_{-2}-x_1-x_{-1} = y_2+y_{-2}-y_1-y_{-1} = 1$.  Consider the monic and root-centered polynomial 
$$F(x) = 24\Psi_4(x+\overline{\alpha})/\binom{v-4}{4} = x^4+p_1 x^3+p_2 x^2 + p_3 x + p_4.$$ 
By Equation (15) in \cite{Peterson}, we have
\begin{align}
p_2 &= -\frac{5}{2}-\frac{6(k-3)(k-4)(v-k-3)(v-k-4)}{(v-6)(v-7)^2},\nonumber\\
\label{p-intermsof-kv}
p_3 &= \frac{-4(k-3)(k-4)(v-k-3)(v-k-4)(v-2k+1)(v-2k-1)}{(v-5)(v-6)(v-7)^3},\\
p_4 &= \frac{9}{16}+\frac{3}{2}\cdot\frac{(k-3)(k-4)(v-k-3)(v-k-4)g(k,v)}{(v-4)(v-5)(v-6)(v-7)^4},\nonumber
\end{align}
where $g(k,v)$ is as in Appendix~\ref{fgkv}. Assuming $x_2+x_{-2}-x_1-x_{-1} = 1$, the roots of $F(x)$ must be $r_1-1/4, r_2+1/4, -r_1-1/4, -r_2+1/4$ where $r_1-1/4 = x_1-\overline{\alpha}$ and $r_2+1/4 = x_2-\overline{\alpha}$ (note that $r_1, r_2 \in \tfrac{1}{4} \mathbb{Z}$). Expanding,
$$x^4+p_2x^2+p_3x+p_4 = x^4+\left(-\tfrac{1}{8}-r_1^2-r_2^2\right)x^2+\left(\tfrac{r_1^2}{2}-\tfrac{r_2^2}{2}\right)x+\left(r_1^2-\tfrac{1}{16}\right)\left(r_2^2-\tfrac{1}{16}\right)\text{,}$$
which yields
\begin{equation}
\label{s4identity}
p_4 = \left(\frac{p_2}{2}+\frac{1}{8}\right)^2-p_3^2.
\end{equation}
Substituting (\ref{p-intermsof-kv}) into (\ref{s4identity}) results in the equation $f(k,v) = 0$.  An easy computer search shows that there are no integer solutions $v$ to $f(k,v)=0$ for $9 \le k \le$ 25,000. 
\qed

By reducing $f(k,v)$ modulo some primes, one may obtain infinite classes of both $k$ and $v$ which admit no soultions.  Some more (very easy) computing is required here.

\begin{cor}
There is no nontrivial tight $8$-$(v,k,\lam)$ design with parameters in any of the following congruence classes:
$$
\begin{array}{l}
k  \equiv  12 \pmod{13} \\
k  \equiv  5,9,11,12 \pmod{17} 
\end{array}
~
\begin{array}{l}
v  \equiv  3 \pmod{7} \\
v  \equiv  3 \pmod{11} \\
v  \equiv  1,3 \pmod{13} \\
v  \equiv  0,1,8,12,13,14 \pmod{17} 
\end{array}
$$
\end{cor}
Despite these strict conditions on hypothetical tight $8$-designs, it remains open whether there are a finite number of nontrivial such designs.

To loosely summarize our work, we have shown that any unknown tight $2s$-design with $s>1$, if it exists, must have
\begin{itemize}
\item
large $s$ and small $\beta$, or
\item
$s=4$ with very large $k$ and $v$ satisfying strict conditions.
\end{itemize}

\section*{Appendices}
\appendix

\section{\large{Hermite polynomials $H_s(x)$ and their zeros, $1 \le s \le 9$}}
\label{hzeros}

\begin{center}
\begin{tabular}{cl}
$s$ & $H_s(x)$ \\
\hline
1 & $x$ \\ 
2 & $x^2-1$ \\
3 & $x^3-3x$ \\
4 & $x^4-6x^2+3$ \\
5 & $x^5-10x^3+15x$ \\
6 & $x^6-15x^4+45x^2-15$ \\
7 & $x^7-21x^5+105x^3-105x$ \\
8 & $x^8-21x^6+210x^4-420x^2+105$ \\
9 & $x^9-36x^7+378x^5-1260x^3+945x$ \\
\hline
\end{tabular}
\end{center}

\noindent
Note: In the following table, the values of $D_i=|H_{s-1}(\xi_i)|$ are rounded down.
$$
\begin{array}{|l|l|l|l|}
\hline
s & i  & \xi_i~\text{for}~H_s(x) & D_i \ge \\ 
\hline
1 & 0 & 0 & 1 \\ 
\hline
2 & 1 & 1 & 1 \\ 
\hline
\multirow{2}{*}{3} & 0 & 0 & 1  \\
		& 1 & \sqrt{3} &  2  \\
\hline
\multirow{2}{*}{4} & 1 & \sqrt{3-\sqrt{6}} = 0.7420 & 1.817  \\
		& 2 & \sqrt{3+\sqrt{6}} = 2.3344 & 5.718  \\
\hline
\multirow{3}{*}{5} & 0 & 0 &  3    \\
    & 1 & \sqrt{5-\sqrt{10}} = 1.3556 & 4.649  \\
		& 2 & \sqrt{5+\sqrt{10}} = 2.8570&  20.64   \\
		\hline
\multirow{3}{*}{6} & 1 & 0.61670659019 & 6.994 \\
    & 2 & 1.88917587775 &  15.02 \\
		& 3 & 3.32425743355 &  88.46 \\
\hline
\end{array}
~
\begin{array}{|l|l|l|l|}
\hline
s & i & \xi_i & D_i \ge \\ 
\hline
\multirow{4}{*}{7} & 0 & 0 & 15  \\
    & 1 & 1.1544 &  20.69 \\
		& 2 & 2.3668 &  57.82 \\
		& 3 & 3.7504 &   433.1  \\
\hline
\multirow{4}{*}{8} & 1 & 0.5391 &  41.09  \\
    & 2 & 1.6365 &  73.30 \\
		& 3 & 2.8025 &  255.7 \\
		& 4 & 4.1445 &  2365 \\
\hline
\multirow{5}{*}{9} & 0 & 0 & 105 \\
    & 1 & 1.0233 &  135.4 \\
		& 2 & 2.0768 &  299.5   \\
		& 3 & 3.2054 &  1267   \\
		& 4 & 4.5127 &  14159   \\
\hline
\end{array}
$$

\section{\large{Constants $B_i$ and $C_i$ obtained from Algorithm~\ref{al} and values of $\beta_0(s)$}}
%%latestmod
\label{C_is}

Notes: For convenience, $B_i$ was chosen independently of $i$ and $C_i$ was taken with $C_i = C_{-i}$.
\begin{center}
\begin{tabular}{|l|l|l|l|l|}
\hline
$s$ & $B_i$ & $C_i$, $i=(0),1,\dots,\fs$ & $\beta_0(s)$ & $\beta_\star(s)$ \\ 
\hline
4 & 10 & $2,14$ & & 19.35 \\%%latestmod
5 & 10 & $1,12,88$ & 33.76& \\ 
6 & 100 & $11,63,558$ & 156.96& \\
7 & 10 & $6,93,458,4649$ & 86.55& \\
8 & 100 & $100,501,2561,30779$ & 106.77& \\
9 & 100 & $9,773,3186,17732,247789$ & 146.37 &\\
\hline
\end{tabular}
\end{center}

\section{\large{The $f(k,v)$ and $g(k,v)$ used in Proposition~\ref{sEquals4}}}
\label{fgkv}
$f(k,v) = -3408102864+1506333312k^{2}+974873344k^{4}-488998144k^{6}+62323584k^{8}-3309568k^{10}+65536k^{12}+9310949028v-1506333312kv-4733985888k^{2}v-1949746688k^{3}v-1015706784k^{4}v+1466994432k^{5}v+511604992k^{6}v-249294336k^{7}v-49810560k^{8}v+16547840k^{9}v+1744896k^{10}v-393216k^{11}v-16384k^{12}v-11097146016v^{2}+4733985888kv^{2}+6922441360k^{2}v^{2}+2031413568k^{3}v^{2}-1428764528k^{4}v^{2}-1534814976k^{5}v^{2}+209662720k^{6}v^{2}+199242240k^{7}v^{2}-21567744k^{8}v^{2}-8724480k^{9}v^{2}+786432k^{10}v^{2}+98304k^{11}v^{2}+7281931941v^{3}-5947568016kv^{3}-4944873072k^{2}v^{3}+412538336k^{3}v^{3}+1856597696k^{4}v^{3}+243542016k^{5}v^{3}-293538048k^{6}v^{3}-13016064k^{7}v^{3}+17194752k^{8}v^{3}-327680k^{9}v^{3}-253952k^{10}v^{3}-2755473732v^{4}+3929166288kv^{4}+1497511456k^{2}v^{4}-1155170432k^{3}v^{4}-582955856k^{4}v^{4}+183266304k^{5}v^{4}+58253568k^{6}v^{4}-16432128k^{7}v^{4}-1102464k^{8}v^{4}+368640k^{9}v^{4}+544096980v^{5}-1459281552kv^{5}+28759472k^{2}v^{5}+469164960k^{3}v^{5}-7038496k^{4}v^{5}-59703552k^{5}v^{5}+6536960k^{6}v^{5}+2050560k^{7}v^{5}-328320k^{8}v^{5}-18769932v^{6}+293023248kv^{6}-127930016k^{2}v^{6}-58917568k^{3}v^{6}+27050224k^{4}v^{6}+1258752k^{5}v^{6}-1642240k^{6}v^{6}+182784k^{7}v^{6}-14780538v^{7}-24513072kv^{7}+27560816k^{2}v^{7}-2875616k^{3}v^{7}-2296192k^{4}v^{7}+698880k^{5}v^{7}-61184k^{6}v^{7}+2961396v^{8}-764688kv^{8}-1582560k^{2}v^{8}+772608k^{3}v^{8}-143664k^{4}v^{8}+10752k^{5}v^{8}-191952v^{9}+203472kv^{9}-52816k^{2}v^{9}+7520k^{3}v^{9}-640k^{4}v^{9}+972v^{10}-2352kv^{10}+336k^{2}v^{10}+45v^{11}$.\\
\\
$g(k,v) = 2k^{4}v-26k^{4}-4k^{3}v^{2}+52k^{3}v+2k^{2}v^{3}-20k^{2}v^{2}-120k^{2}v+258k^{2}-6kv^{3}+120kv^{2}-258kv+v^{4}-23v^{3}+123v^{2}-433v+764$.

\section*{Acknowledgement}
The authors would like to thank Jane Wodlinger for helpful discussions.

\end{document}